\documentclass[12pt]{amsart}
\usepackage{amssymb, latexsym, amsmath, amscd, oldlfont}
\usepackage[dvips]{graphicx}
\usepackage{psfrag}

\advance\hoffset by -0.5 truecm

\def\V{\widehat{V}}
\def\U{\widehat{U}}

\def\H{{\cal H}}
\def\HH{{\widehat{H}}}
\def\K{\widehat{K}}

\def\ppi^n{\widetilde{\pi}}
\def\al{{\alpha}}

\def\ep{{\epsilon}}
\def\ga{{\gamma}}

\def\ov{\overline}
\def\om{{\omega}}
\def\Om{{\Omega}}

\def\te{{\theta}}

\def\lg{\langle}
\def\rg{\rangle}

\def\vr{\varphi}

\def\pr{\prime}

\def\ppi{\widetilde{\pi}}

\def\lra{\longrightarrow}

\def\sm{\setminus}
\def\Sig{\Sigma}
\def\hra{\hookrightarrow}
\def\Hp{{H\circ\psi}}
\def\gga{{\widetilde{\gamma}}}
\def\FF{{\widetilde{F}}}
\def\FFF{{\widehat{F}}}
\def\JJ{{\widetilde{J}}}
\def\gga{{\widehat{\gamma}}}
\def\HZ{\operatorname{HZ}}
\def\cHZ{{\widetilde{c}}_{\HZ}}
\def\Ht{{\widetilde{{\cal H}}}}

\def\C{{\mathbb C}}
\def\R{{\mathbb R}}
\def\N{{\mathbb N}}
\def\Z{{\mathbb Z}}

\newtheorem{theorem}{Theorem}[section]
\newtheorem{proposition}{Proposition}[section]
\newtheorem{lemma}{Lemma}[section]
\newtheorem{corollary}{Corollary}[section]
\theoremstyle{definition}
\newtheorem{definition}{Definition}[section]

\def\thebibliography#1{\section*{References}\list
 {[\arabic{enumi}]}{\settowidth\labelwidth{[#1]}\leftmargin\labelwidth
 \advance\leftmargin\labelsep
 \usecounter{enumi}}
 \def\newblock{\hskip .11em plus .33em minus .07em}
 \sloppy\clubpenalty4000\widowpenalty4000
 \sfcode`\.=1000\relax}

\newenvironment{thmintro}{\par\medskip\noindent{\bf 
Main Theorem.}\begingroup\em}{\endgroup\hfill\par\medskip}

\newenvironment{remark}{\par\medskip\noindent{\bf 
Remark.}}{\hfill\par\medskip}
\newenvironment{remarks}{\par\medskip\noindent{\bf 
Remarks.}}{\hfill\par\medskip}

\author{Leonardo Macarini}
\title[Hofer-Zehnder semicapacity of cotangent bundles]{Hofer-Zehnder semicapacity of cotangent bundles and symplectic submanifolds}
\address{Instituto de Matem\'atica Pura e Aplicada - IMPA\\
         Estrada Dona Castorina, 110 - Jardim Bot\^anico\\
         22460-320 Rio de Janeiro RJ\\
         Brasil}
\email{leonardo@impa.br}         
\date{current version: June 2003}
\thanks{This work was partially supported by CNPq-PROFIX, Brazil.}

\begin{document}

\begin{abstract}
We introduce the concept of Hofer-Zehnder $G$-semicapacity (or $G$-sensitive Hofer-Zehnder capacity) and prove that given a geometrically bounded symplectic manifold $(M,\om)$ and an open subset $N \subset M$ endowed with a Hamiltonian free circle action $\vr$ then $N$ has bounded Hofer-Zehnder $G_\vr$-semicapacity, where $G_\vr \subset \pi_1(N)$ is the subgroup generated by the homotopy class of the orbits of $\vr$. In particular, $N$ has bounded Hofer-Zehnder capacity.

We give two types of applications of the main result. Firstly, we prove that the cotangent bundle of a compact manifold endowed with a free circle action has bounded Hofer-Zehnder capacity. In particular, the cotangent bundle $T^*G$ of any compact Lie group $G$ has bounded Hofer-Zehnder capacity. Secondly, we consider Hamiltonian circle actions given by symplectic submanifolds. For instance, we prove the following generalization of a recent result of Ginzburg-G\"urel \cite{GG}: almost all levels of a function defined on a neighborhood of a closed symplectic submanifold $S$ in a geometrically bounded symplectic manifold carry contractible periodic orbits of the Hamiltonian flow, provided that the function is constant on $S$.
\end{abstract}

\maketitle

\section{Introduction}

We will consider here the problem of the existence of periodic orbits on prescribed energy levels of Hamiltonian systems. An approach to this problem was the introduction of certain symplectic invariants related to symplectic rigidity phenomena. More precisely, H. Hofer and E. Zehnder introduced in \cite{HZ} the following definition:

\begin{definition}
\label{defcHZ}
Given a symplectic manifold $(M,\om)$, define the Hofer-Zehnder capacity of $M$
by
$$ c_{\HZ}(M,\om) = \sup_{H \in \H_a(M,\om)}\ \max H, $$
where $\H_a(M,\om)$ is the set of {\it admissible Hamiltonians} $H$ defined on
$M$, that is,
\begin{itemize}
\item $H \in \H(M) \subset C^\infty(M,\R)$, where $\H(M)$ is the set
of {\it pre-admissible Hamiltonians} defined on $M$, that is, $H$ satisfies the
following properties: $0 \leq H \leq \|H\|:=\max H - \min H$, there exist an open set 
$V \subset M$ such that $H|_V \equiv 0$ and a compact set 
$K \subset M\setminus\partial M$ satisfying $H|_{M\setminus K}\equiv \|H\|$;
\item every non-constant periodic orbit of $X_H$ has period greater than 1.
\end{itemize}
\end{definition}

We will also consider here a modified Hofer-Zehnder capacity introduced by D. McDuff and J. Slimowitz \cite{MSl} sensitive to the presence of overtwisted critical points. It is given by
$$ \cHZ(M,\om) = \sup_{H \in \Ht_a(M,\om)}\ \max H, $$
where $H \in \Ht_a(M,\om)$ if $H$ is admissible and has no overtwisted critical points, that is,
\begin{itemize}
\label{hyp}
\item the linearized flow at the critical points of $H$ has no non-trivial periodic orbits of period less than 1.
\end{itemize}

\begin{remark}
It is an open question whether the capacities $c_{HZ}$ and $\cHZ$ coincide.
\end{remark}

It is easy to prove that if a symplectic manifold $(M,\om)$ has {\it bounded Hofer-Zehnder capacity}, that is, if $c_{\HZ}(U,\om) < \infty$ for every open subset $U \subset M$ with compact closure, then given any Hamiltonian $H: M \to \R$ with compact energy levels, there exists a dense subset $\Sigma \subset H(M)$ such that for every $e \in \Sigma$ the energy hypersurface $H^{-1}(e)$ has a periodic solution. Actually, if there exists a neighborhood of an energy hypersurface $H^{-1}(e_0)$ with finite Hofer-Zehnder capacity then there exists $\ep>0$ and a dense subset $G \subset (e_0-\ep,e_0+\ep)$ such that $H^{-1}(e)$ has periodic orbits for every $e \in G$. This result was improved in \cite{MS}, where it is showed that there are periodic orbits on $H^{-1}(e)$ for almost all $e \in (e_0-\ep,e_0+\ep)$ with respect to the Lebesgue measure. The same holds for $\cHZ(M,\om)$.

Let us consider here a refinement of the original Hofer-Zehnder capacity by
considering periodic orbits whose homotopy class is contained in a given
subgroup $G$ of $\pi_1(M)$. More precisely, we define the Hofer-Zehnder 
semicapacity as follows:

\begin{definition}
Given a symplectic manifold $(M,\om)$ and a subgroup $G \subset \pi_1(M)$
define the {\it Hofer-Zehnder G-semicapacity} (or {\it $G$-sensitive Hofer-Zehnder capacity}) of $M$ by
$$ c_{\HZ}^G(M,\om) = \sup_{H \in \H_a^G(M,\om)}\ \max H, $$
where $\H_a^G(M,\om)$ is the set of {\it $G$-admissible Hamiltonians} defined
on $M$, that is,
\begin{itemize}
\item $H \in \H(M)$, that is, $H$ is pre-admissible (see the definition
\ref{defcHZ});
\item every nonconstant periodic orbit of $X_H$ whose homotopy class belongs to
$G$ has period greater than 1.
\end{itemize}
We define $\cHZ^G(M,\om)$ analogously considering $G$-admissible Hamiltonians without overtwisted critical points.
\end{definition}

\begin{remark}
We call it a {\it semicapacity} because given a symplectic embedding  $\psi:
(N,\tau) \to (M,\om)$ such that $\dim N = \dim M$ it cannot be expected, in
general, that $c_{\HZ}^G(N,\tau) \leq c_{\HZ}^{\psi_*G}(M,\om)$, where
$\psi_*: \pi_1(N) \to \pi_1(M)$ is the induced homomorphism on the fundamental
group. However, we can state the following {\it weak monotonicity property}
(for a proof see \cite{Mac}): given a symplectic embedding
$\psi: (N,\tau)\to (M,\om)$ such that $\dim N = \dim M$ then
$$ c_{\HZ}^{\psi_*^{-1}H}(N,\tau) \leq c_{\HZ}^H(M,\om),$$
for every subgroup $H \subset \pi_1(M)$. In particular, if $\psi_*$ is injective we have that  $c_{\HZ}^G(N,\tau)
\leq c_{\HZ}^{\psi_*G}(M,\om)$.
\end{remark}

Note that obviously,
$$ c_{\HZ}(M,\om) = c_{\HZ}^{\pi_1(M)}(M,\om) \leq c_{\HZ}^G(M,\om), $$
for every subgroup $G \subset \pi_1(M)$. Moreover, it can be show that, like the Hofer-Zehnder capacity, if
the Hofer-Zehnder $G$-semicapacity is bounded then there are periodic orbits
with homotopy class in $G$ on almost all energy levels for every proper
Hamiltonian \cite{MS}.

\begin{definition}
A symplectic manifold $(M,\om)$ is called {\it geometrically bounded} if there exists an almost complex structure $J$ on $M$ and a Riemannian metric $g$ such that
\begin{itemize}
\item $J$ is uniformly $\om$-tame, that is, there exist positive constants $c_1$ and $c_2$ such that
$$ \om(v,Jv) \geq c_1\|v\|^2 \text{ and } |\om(v,w)| \leq c_2\|v\|\|w\| $$
for all tangent vectors $v$ and $w$ to $M$;
\item the sectional curvature of $g$ is bounded from above and the injectivity radius is bounded away from zero.
\end{itemize}
\end{definition}

Closed symplectic manifolds are clearly geometrically bounded; a product of two geometrically bounded symplectic manifolds is also such a manifold. It was showed in \cite{CGK,Lu1} that the cotangent bundle of a compact manifold endowed with any twisted symplectic form is geometrically bounded. A {\it twisted symplectic structure} on a cotangent bundle $T^*M$ is a symplectic form given by $\om_0 + \pi^*\Om$, where $\om_0$ is the canonical symplectic form, $\pi: T^*M \to M$ is the canonical projection and $\Om$ is a closed 2-form on $M$.

\begin{definition}
Given a symplectic manifold $(M,\om)$ define its {\it index of rationality} by
$$ m(M,\om) = \inf\bigg\{\int_{[u]}\om;\ \ [u] \in \pi_2(M) \text{ satisfies } \int_{[u]}\om > 0\bigg\}. $$
If the above set is empty, we define $m(M,\om) = \infty$. If $m(M,\om)>0$ we call $(M,\om)$ a {\it rational symplectic manifold}.
\end{definition}

The main result here states that the existence of a free Hamiltonian circle action implies the boundedness of the Hofer-Zehnder semicapacity. It removes a technical hypothesis in a previous theorem in \cite{Mac}. 

\begin{thmintro}
\label{thmA} 
Let $(M,\om)$ be a geometrically bounded symplectic manifold and $N \subset M$ an open subset that admits a free Hamiltonian circle action $\vr$. Then, given any open subset $U \overset{i}{\hra} N$ with compact closure,
$$ \cHZ^{i_*^{-1}G_\vr}(U,\om) < \infty, $$
where $i_*: \pi_1(U) \to \pi_1(N)$ is the induced homomorphism on the
fundamental group and $G_\vr \subset \pi_1(N)$ is the subgroup generated by the
orbits of the circle action. Moreover, if $M$ is rational, then $c^{i_*^{-1}G_\vr}_{\HZ}(U,\om) < \infty$.
\end{thmintro}

The main idea in the proof of the Main Theorem is to relate the Hofer-Zehnder capacity of a symplectic manifold endowed with a free Hamiltonian action with the Hofer-Zehnder capacity of its symplectic reduction. More precisely, let $P$ be a symplectic manifold with a free Hamiltonian $G$-action and bounded Hofer-Zehnder capacity. When its reduced symplectic manifold $M$ also has bounded Hofer-Zehnder capacity? Of course, it does not hold in general, since every symplectic manifold $M$ can be obtained as a reduction of $M \times T^*S^1$ with respect to the obvious circle action on $T^*S^1$ and this symplectic manifold has bounded Hofer-Zehnder capacity by the Main Theorem. On the other hand, the example of Zehnder \cite{Ze} shows that there exist compact symplectic manifolds with infinite Hofer-Zehnder capacity.

Let us consider the easiest case, namely when $G=S^1$. In what follows, we will always denote the symplectic gradient of a Hamiltonian $H$ by $X_H$. Given a pre-admissible Hamiltonian $H$ on $M$ we can construct a pre-admissible Hamiltonian $\HH$ on $P$ whose reduced dynamics coincides with the dynamics of $H$. The problem is to give conditions to ensure that the reduced periodic orbit on $M$ is not trivial, that is, that the periodic orbit of $X_\HH$ is not tangent to the orbits of the group action. An approach to solve this problem is to use the fact that besides the existence of a non-trivial periodic orbit, the Hofer-Zehnder capacity gives an upper bound $T$ for the period of the orbit. Thus, we can try to construct a circle action with sufficiently great period in such a way that the Hofer-Zehnder capacity of the support of $X_\HH$ is sufficiently small (compared to the period of the action) and that the component of $X_\HH$ tangent to the circle action remains bounded. It would ensure that if a periodic orbit of $X_\HH$ is tangent to the orbits of the action then its period is necessarily greater than $T$.

When $M$ has a free Hamiltonian circle action we construct such a circle action on $M \times T^*S^1$ but, instead of $M$, its symplectic reduction is given by a finite quotient $M/\Z_n$ for $n$ sufficiently great (Theorem \ref{thm1}). We then get the boundedness of the Hofer-Zehnder semicapacity of $M/\Z_n$ (resp. $\cHZ^G$) by Hofer-Viterbo theorem \cite{HV} (resp. McDuff-Slimowitz theorem \cite{MSl}) generalized for geometrically bounded symplectic manifolds by G. Lu \cite{Lu1,Lu2}. It was essentially achieved in \cite{Mac}. However, we need a theorem to relate the Hofer-Zehnder capacity of $M$ and its finite symplectic quotient $M/\Z_n$. This is the content of Theorem \ref{thm2} proved in Section \ref{proofthm2}.

\subsection{Applications to cotangent bundles}
Now, let us consider some applications of the Main Theorem. Considering the lifted action to the cotangent bundle, we have the following corollary:

\begin{corollary}
Let $M$ be a manifold with a free circle action $\vr: S^1 \times M \to M$. Then $T^*M$ endowed with the canonical symplectic form has bounded Hofer-Zehnder $G_\vr$-semicapacity, where $G_\vr \subset \pi_1(M) = \pi_1(T^*M)$ is the subgroup generated by the orbits of the action.
\end{corollary}

Now, let $G$ be a compact Lie group and consider the action of the maximal torus in $G$. It gives the following immediate corollary:

\begin{corollary}
Let $G$ be a compact Lie group. Then $T^*G$ has bounded Hofer-Zehnder capacity.
\end{corollary}

\subsection{Applications to symplectic submanifolds}
The next corollary is a generalization of a recent result of V. Ginzburg and B. G\"urel \cite{GG} for weakly exact symplectic manifolds. Their proof is completely different and is based on symplectic homology. Moreover, they consider only energy hypersurfaces that bound a compact domain containing the symplectic submanifold $S$ while our result holds for any energy hypersurface that does not intersect $S$.

\begin{corollary}
Let $(M,\om)$ be a geometrically bounded symplectic manifold and $S \subset M$ a closed symplectic submanifold. There exists a neighborhood $U$ of $S$ such that if $H$ is a proper Hamiltonian on $U$ constant on $S$ then $H$ has periodic orbits with contractible projection on $S$ on almost all energy levels.
\end{corollary}

\begin{remark}
The only condition on $U$ is that $U\setminus S$ admits a free Hamiltonian circle action. Its existence follows by the symplectic neighborhood theorem \cite{MSa}.
\end{remark}

As an immediate corollary we have the following result about the existence of periodic orbits for magnetic flows (for an introduction to magnetic flows see \cite{CMP,Mac}):

\begin{corollary}
Let $(M,\Om)$ be a closed symplectic manifold and $H: T^*M \to \R$ be the standard kinetic energy Hamiltonian for a Riemannian metric on $M$. Then the Hamiltonian flow of $H$ with respect to the twisted symplectic form $\om_0 + \pi^*\Om$ has periodic orbits with contractible projection on $M$ on almost all sufficiently low energy levels.
\end{corollary}

The next application is based on the following result of P. Biran \cite{Bir} which enable us to represent a K\"ahler manifold as a disjoint union of two basic components whose symplectic nature is very standard:

\begin{theorem}[P. Biran \cite{Bir}]
\label{Bir}
Let $(M^{2n},\Om)$ be a closed K\"ahler manifold with $[\Om] \in
H^2(M,\Z)$ and $\Sigma \subset M$ a complex hypersurface whose homology class
$[\Sigma] \in H_{2n-2}(M)$ is the Poincar\'e dual to $k[\Om]$ for some $k \in
\N$. Then, there exists an isotropic CW-complex $\Delta \subset (M,\Om)$ whose
complement - the open dense subset $(M\sm\Delta,\Om)$ - is symplectomorphic to
a standard symplectic disc bundle $(E_0,\frac{1}{k}\om_{\text{can}})$ modeled on the
normal bundle $N_\Sig$ of $\Sig$ in $M$ and whose fibers have area $1/k$.
\end{theorem}

\begin{remark}
The CW-complex $\Delta$ above is given by the union of the stable manifolds of the gradient flow of a plurisubharmonic function defined on $M\sm\Sig$. Thus, it can be explicitly computed in many examples \cite{Bir}.
\end{remark}

The symplectic form $\om_{\text{can}}$ is given by
$$ \om_{\text{can}} = k\pi^*(\Om|_\Sigma) + d(r^2\al), $$
where $\pi: E_0 \to \Sigma$ is the bundle projection, $r$ is the radial
coordinate using a Hermitian metric $\|\cdot\|$ and $\al$ is a connection form
on $E$ such that $d\al = -k\pi^*(\Om|_\Sigma)$. The form
$\frac{1}{k}\om_{\text{can}}$ is uniquely characterized by the requirements
that its restriction to the zero section $\Sig$ equals $\Om|_\Sig$, the fibers
of $\pi: E_0 \to \Sigma$ are symplectic and have area $1/k$ and
$\om_{\text{can}}$ is invariant under the obvious circle action along the
fibers. It is called standard because the symplectic type of
$(E_0,\om_{\text{can}})$ depends only on the symplectic type of
$(\Sig,\Om|_\Sig)$ and the topological type of the normal bundle $N_\Sig$
\cite{Bir}.

Let us recall that the pair $(M,\Sig)$ is called {\it subcritical} \cite{BC1,BC2} if $M\sm\Sig$ is
a subcritical Stein manifold, that is, if there exists a plurisubharmonic Morse
function $\vr$ on $M\sm\Sig$ such that $\text{index}_p(\vr) < \dim_\C M$ for
every critical point $p$ of $\vr$. It is equivalent to the condition that
the dimension of $\Delta$ (that is, the maximal dimension of the cells of
$\Delta$) is strictly less than $n$.

Since $(E_0\setminus\Sig,\om_{\text{can}})$ has an obvious Hamiltonian free circle action, we have the following immediate corollary of the Main Theorem:

\begin{corollary}
According to the notation of the Theorem \ref{Bir}, we have that
$$ c^G_{\HZ}(E_0\sm\Sig,\Om) < \infty, $$
where $G \subset \pi_1(E_0\sm\Sigma)$ is the subgroup generated by the orbits of
the obvious $S^1$-action on $E_0\sm\Sig$. In particular, the periodic orbits are contractible in $M$. Moreover, if $(M,\Sig)$ is subcritical, then 
$$ c^{i_*G}_{\HZ}(M\sm\Sigma,\Om) < \infty, $$
where $E_0\sm\Sig \overset{i}{\hra} M\sm\Sig$.
\end{corollary}

\begin{remarks}
\begin{itemize}
\item The result for subcritical manifolds follows by the fact that if $\dim \Delta < \frac{1}{2}\dim_\R M$ then there exists a Hamiltonian isotopy $k_t: (M\sm\Sig,\Om) \to (M\sm\Sig,\Om)$ compactly supported in an arbitrarily small neighborhood of $\Delta$ such that $k_1(\Delta) \cap \Delta = \emptyset$. For the details see \cite{Bir,Mac}.
 
\item It removes the hypothesis $\Om|_{\pi_2(M)}=0$ in Theorem 1.4 of \cite{Mac}.
\end{itemize}
\end{remarks}


\medskip
\noindent
{\bf Acknowledgements.} I am very grateful to Felix Schlenk for very useful comments.

\section{Proof of the Main Theorem}

The proof is based on the following theorems whose proofs are given in the Sections \ref{proofthm1} and \ref{proofthm2} respectively:

\begin{theorem}
\label{thm1}
Let $(M,\om)$ be a geometrically bounded symplectic manifold and $N \subset M$ an open subset that admits a free Hamiltonian circle action $\vr$ (with period equal to 1) generated by the Hamiltonian $H_1$. Consider on $N$ the $\Z_n$-action given by $\vr_{1/n}$. Given an open subset $U \subset N$ with compact closure, there exists a positive integer $n_0=n_0(\|H_1|_U\|) < \infty$ such that given $n\geq n_0$ and a pre-admissible Hamiltonian $H$ satisfying
\begin{itemize}
\item $H$ has no overtwisted critical points;
\item $H$ is $\Z_n$-invariant;
\item the support of $X_H$ is contained in $U$;
\item $ \|H\| > 2\pi(\|H_1|_U\|+\sqrt{n}), $
\end{itemize}
then there exists a non-trivial periodic orbit of $X_H$ with period less than 1 and homotopy class in $G_\vr$. Moreover, if $M$ is rational, the hypothesis that $H$ has no overtwisted critical points is not necessary, but $n_0$ also depends on $m(M,\om)$.
\end{theorem}

\begin{theorem}
\label{thm2}
Let $n = 2^m$ for some positive integer $m$ and $M$ be an open  symplectic manifold endowed with a free symplectic $\Z_n$-action, that is, a symplectomorphism $\psi: M \to M$ such that $\psi^n=Id$ and $\psi^i(x)\neq x$ for every $x \in M$ and $1\leq i < n$. Then, given an admissible Hamiltonian $H$ on $M$ (resp. $H \in \Ht_a(M,\om)$), there exists an admissible Hamiltonian $H^\prime$ (resp. $H^\prime \in \Ht_a(M,\om)$) such that
\begin{itemize}
\item $H^\prime$ is $\Z_n$-invariant;
\item $\text{supp }X_{H^\prime} = \text{supp }X_{\sum_{i=1}^{n}H\circ\psi^i}$;
\item $\|H^\prime\| \geq (1/n)\|H\|$.
\end{itemize}
Moreover, if $\psi$ is isotopic to the identity and $H$ is $G$-admissible (resp. $H \in \Ht_a^G(M,\om)$) then $H^\prime$ is also $G$-admissible (resp. $H^\prime \in \Ht_a^G(M,\om)$), for any subgroup $G \subset \pi_1(M)$.
\end{theorem}

\begin{proof}[Proof of the Main Theorem: ]
Suppose, without loss of generality, that the orbits of $\vr$ have period 1. Let $n = 2^m$ be such that $n \geq n_0$, where $n_0=n_0(\|H_1|_U\|,m(M,\om))$ is the integer given by Theorem \ref{thm1}. Let $\psi(x) = \vr_{1/n}(x)$ be the symplectic $\Z_n$-action induced by $\vr$ and $H$ a $G_\vr$-admissible Hamiltonian on $U$ without overtwisted critical points. By Theorem \ref{thm2}, there exists a $\Z_n$-invariant $G_\vr$-admissible Hamiltonian $H^\prime$ without overtwisted critical points such that
$$\|H^\prime\| \geq (1/n)\|H\|. $$
Moreover, since $\text{supp }X_{H^\prime} = \text{supp }X_{\sum_{i=1}^{n}H\circ\psi^i}$, we have that
$$ \|H_1|_{\text{supp }X_{H^\prime}}\| = \|H_1|_{\text{supp }X_H}\| \leq \|H_1|_U\|. $$
Thus, we conclude that if
$$ \|H\| > 2\pi n(\|H_1|_U\|+\sqrt{n}), $$
then $\|H^\prime\| > 2\pi(\|H_1|_U\|+\sqrt{n})$ and consequently, by Theorem \ref{thm1}, there exists a non-trivial periodic orbit of $X_{H^\prime}$ with homotopy class in $G_\vr$ and period less than 1, contradicting the fact that $H^\prime$ is $G_\vr$-admissible. So, we have that
$$ \cHZ^{i_*^{-1}G_\vr}(U,\om) \leq 2\pi n(\|H_1|_U\|+\sqrt{n}) < \infty. $$
When $M$ is rational, the hypothesis on the critical points is not necessary and hence
$$ c_{\HZ}^{i_*^{-1}G_\vr}(U,\om) \leq 2\pi n(\|H_1|_U\|+\sqrt{n}) < \infty. $$
\end{proof}

\section{Proof of Theorem \ref{thm1}}
\label{proofthm1}

If $H$ has no overtwisted critical points, let $n_0 = n_0(\|H_1|_U\|)$ be given by
$$ n_0 = \inf\{n \in \N;\ \sqrt{n} > (4\pi/n)(\|H_1|_U\| + \sqrt{n})\}. $$
If $H$ has overtwisted critical points and $M$ is rational, let $n_0 = n_0(\|H_1|_U\|,m(M,\om))$ be the infimum of $n \in \N$ such that $\sqrt{n} > (4\pi/n)(\|H_1|_U\| + \sqrt{n})$ and
$$ m(M,\om) \geq (2\pi/n)(\|H_1|_U\|+\sqrt{n}). $$
We will divide the proof in three steps:

\subsection{Construction of a circle action $\rho$ on $N \times T^*S^1$ and its symplectic reduction. }
\label{step1}
Let $P = N \times T^*S^1 \subset M \times T^*S^1$, where $M \times T^*S^1$ is endowed with the symplectic form $\om_P = \om \oplus \om^{S^1}_0$ and $\om^{S^1}_0$ is the canonical symplectic form on $T^*S^1$. Note that $T^*S^1$ has a natural Hamiltonian circle action (with period equal to 1) given by the Hamiltonian $\Phi(\theta,\mu)=\mu$, where $\theta$ is the angle coordinate on $S^1$. Consider the free circle action $\rho: P \times S^1 \to P$ on $P$ given by the Hamiltonian
$$ J(x,\theta,\mu) = (1/\sqrt{n})H_1(x) + \sqrt{n}\Phi(\theta,\mu). $$

We will show that the reduced symplectic manifold $J^{-1}(\mu)/S^1$ with respect to this circle
action is given by $(N/\Z_n,\omega_n)$ for every $\mu \in \R$, where the action of $\Z_n \subset S^1$
is that one induced by $\vr$ and the symplectic form $\om_n$ on $N/\Z_n$ is the
unique form such that $\tau_n^*\om_n=\om$, where $\tau_n: N \to N/\Z_n$ is the
quotient projection (such that $\tau_n: (N,\om) \to (N/\Z_n,\om_n)$ defines a
finite symplectic covering).

In fact, $P/S^1$ can be given by the quotient $(P/\Z_n)/(S^1/\Z_n)$. But, note that the $\Z_n$-action induced by $\rho$ coincides with the $\Z_n$-action induced by $\vr$. Actually, since $\rho$ has period $\sqrt{n}$, the induced $\Z_n$-action is given by $\rho_{\sqrt{n}/n}=\rho_{1/\sqrt{n}}$. On the other hand, $\sqrt{n}\Phi$ also has period $1/\sqrt{n}$ and, since $H_1$ and $\Phi$ commute, we have that $\rho_{1/\sqrt{n}} = \vr_{1/n}$ (see the figure below). Thus,
$$ P/\Z_n = N/\Z_n \times T^*S^1 =: P_n. $$

\begin{figure}[h]
\label{actionJ}
\begin{center}
\psfrag{XH1}{$\frac{1}{\sqrt{n}}X_{H_1}$}
\psfrag{XPhi}{$\sqrt{n}X_\Phi$}
\includegraphics[width=2in]{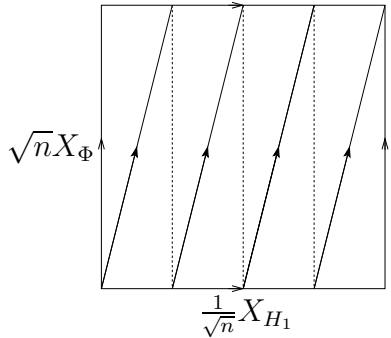}
\caption{\label{actionrho} The action of $X_J$ for $n=4$.}
\end{center}
\end{figure}

Consequently, we need to show that the Marsden-Weinstein reduction $\JJ^{-1}(\mu)/S^1$ with respect to the induced circle action on $P_n$ is given by $(N/\Z_n,\omega_n)$, where $\JJ$ is given by
$$ \JJ(x,\theta,\mu) = (1/\sqrt{n})H_1(\tau_n^{-1}(x)) +  \sqrt{n}\Phi(\theta,\mu). $$
But, since the Hamiltonians $(1/\sqrt{n})H_1(\tau_n^{-1}(x))$ and $\sqrt{n}\Phi(\theta,\mu)$ generate a circle action with the same period $1/\sqrt{n}$, the circle action induced by $\JJ$ is given by the diagonal action. Now, we will use the following proposition:

\begin{proposition}
Let $(Q,\tau)$ be a symplectic manifold endowed with a Hamiltonian circle action given by the Hamiltonian $H_1$ and consider on $(Q \times T^*S^1,\tau\oplus\om_0^{S^1})$ the diagonal circle action generated by $\JJ:= H_1 + \Phi$ (such that $H_1$ and $\Phi$ have the same period). Then the Marsden-Weinstein reduced symplectic manifold $(\JJ^{-1}(\mu)/S^1,\sigma_\mu)$ is $(Q,\tau)$ for every $\mu \in \R$.
\end{proposition}

\begin{proof}
Define $W:= Q \times T^*S^1$ and let $\om_W:=\tau\oplus\om_0^{S^1}$ and $\psi: W \to W$ be the diffeomorphism given by
$$ \psi(x,\theta,\mu) = (\vr_\theta(x),\theta,\mu - H_1(x)), $$
where $\vr: Q \times S^1 \to Q$ is the Hamiltonian flow of $H_1$. The following lemma is a straightforward computation.

\begin{lemma}
The diffeomorphism $\psi$ is a symplectomorphism with respect to $\om_W$.
Moreover, $\JJ \circ \psi = \Phi$.
\end{lemma}

\begin{proof}
It is clear that $\psi^*\JJ = \Phi$. In effect,
\begin{align*}
\JJ(\psi(x,\te,\mu)) & = H_1(\vr_\te(x)) + \mu - H_1(x) \\
& = H_1(x) - H_1(x) + \mu \\
& = \mu.
\end{align*}

To show that $\psi^*\om_W = \om_W$, it is more convenient to write $\psi$ as 
$$ \psi(z) = \vr_{\theta(z)}(z) - H_1(z)Y, $$
where $\theta: W \to S^1$ is the projection onto the circle and
$Y$ is the vector field tangent to the fibers of $T^*S^1$ such that
$\Phi(Y)=1$. Consequently,
$$ d\psi(z)\xi = (d\vr)_{\te(z)}(z)\xi + \pi_2^*d\te(\xi)X_{H_1}(\vr_{\te(z)}(z)) -
dH_1(z)\xi Y(\vr_{\te(z)}(z)), $$
where $X_{H_1}$ is the Hamiltonian vector field generated by $H_1$ and
$\pi_2^*d\te$ is the pullback of the angle form $d\te$ on $S^1$ to $W$. Thus,
\begin{align*}
& (\psi^*\om_W)_z(\xi,\eta) = (\om_W)_{\psi(z)}(d\psi(z)\xi,d\psi(z)\eta) \\
& = \om_W((d\vr)_{\te(z)}(z)\xi,(d\vr)_{\te(z)}(z)\eta) 
+ \om_W((d\vr)_{\te(z)}(z)\xi,\pi_2^*d\te(\eta)X_{H_1}(\vr_{\te(z)}(z))) \\
& - \om_W((d\vr)_{\te(z)}(z)\xi,dH_1(z)\eta Y(\vr_{\te(z)}(z)))
+ \om_W(\pi_2^*d\te(\xi)X_{H_1}(\vr_{\te(z)}(z)),(d\vr)_{\te(z)}(z)\eta) \\
& + \om_W(\pi_2^*d\te(\xi)X_{H_1}(\vr_{\te(z)}(z)),\pi_2^*d\te(\eta)X_{H_1}(\vr_{\te(z)}(z)))\\
& -\om_W(\pi_2^*d\te(\xi)X_{H_1}(\vr_{\te(z)}(z)),dH_1(z)\eta Y(\vr_{\te(z)}(z)))\\
& - \om_W(dH_1(z)\xi Y(\vr_{\te(z)}(z)),(d\vr)_{\te(z)}(z)\eta)
-\om_W(dH_1(z)\xi Y(\vr_{\te(z)}(z)),\pi_2^*d\te(\eta)X_{H_1}(\vr_{\te(z)}(z)))\\
& + \om_W(dH_1(z)\xi Y(\vr_{\te(z)}(z)),dH_1(z)\eta Y(\vr_{\te(z)}(z))).
\end{align*}
But note that
$$\om_W(\pi_2^*d\te(\xi)X_{H_1}(\vr_{\te(z)}(z)),\pi_2^*d\te(\eta)X_{H_1}(\vr_{\te(z)}(z)))=0$$
and
$$\om_W(dH_1(z)\xi Y(\vr_{\te(z)}(z)),dH_1(z)\eta Y(\vr_{\te(z)}(z)))=0,$$
since the vectors are colinear. On the other hand,
$$\om_W(\pi_2^*d\te(\xi)X_{H_1}(\vr_{\te(z)}(z)),dH_1(z)\eta Y(\vr_{\te(z)}(z)))=0$$
and
$$\om_W(dH_1(z)\xi Y(\vr_{\te(z)}(z)),\pi_2^*d\te(\eta)X_{H_1}(\vr_{\te(z)}(z)))=0,$$
because the vectors are orthogonal with respect to the product decomposition $P
= M \times T^*S^1$ and by the definition of $\om_W$.

We have then that,
\begin{align*}
& (\psi^*\om_W)_z(\xi,\eta)
= \om_W((d\vr)_{\te(z)}(z)\xi,(d\vr)_{\te(z)}(z)\eta) 
+ \om_W((d\vr)_{\te(z)}(z)\xi,\pi_2^*d\te(\eta)X_{H_1}(\vr_{\te(z)}(z)))\\
&- \om_W((d\vr)_{\te(z)}(z)\xi,dH_1(z)\eta Y(\vr_{\te(z)}(z)))
+ \om_W(\pi_2^*d\te(\xi)X_{H_1}(\vr_{\te(z)}(z)),(d\vr)_{\te(z)}(z)\eta)\\
&-\om_W(dH_1(z)\xi Y(\vr_{\te(z)}(z)),(d\vr)_{\te(z)}(z)\eta)=\\
& = (\om_W)_z(\xi,\eta) +
\pi_2^*d\te(\eta)(i_{X_{H_1}(\vr_{\te(z)}(z))}\om)(d\vr_{\te(z)}(z)\xi)
- dH_1(z)\eta(i_{Y(\vr_{\te(z)}(z))}\om)(d\vr_{\te(z)}(z)\xi)\\
& + \pi_2^*d\te(\xi)(i_{X_{H_1}(\vr_{\te(z)}(z))}\om)(d\vr_{\te(z)}(z)\eta)
 - dH_1(z)\xi(i_{Y(\vr_{\te(z)}(z))}\om)(d\vr_{\te(z)}(z)\eta),
\end{align*}
where the last equality follows from the fact that $\vr_t^*\om_W=\om_W$. Now,
note that
$$ \pi_2^*d\te(\eta)(i_{X_{H_1}(\vr_{\te(z)}(z))}\om)(d\vr_{\te(z)}(z)\xi)
= \pi_2^*d\te(\eta)dH_1(z)\xi, $$
because $i_{X_{H_1}}\om_W = dH_1$ and $\vr^*dH_1 = dH_1$. On the other hand, we
have that
$$ dH_1(z)\eta(i_{Y(\vr_{\te(z)}(z))}\om)(d\vr_{\te(z)}(z)\xi)
= dH_1(z)\eta \pi_2^*d\te(\xi), $$
because $i_{Y}\om_W = \pi_2^*d\te$ and $\vr^*\pi_2^*d\te = \pi_2^*d\te$ (note
that $\pi_2^*d\te$ here is the pullback of the angle form $d\te$ on $S^1$ to $P$).
Consequently,
\begin{align*}
(\psi^*\om_W)_z(\xi,\eta) & = (\om_W)z(\xi,\eta) + \\
& + \pi_2^*d\te(\eta)dH_1(z)\xi
- dH_1(z)\eta \pi_2^*d\te(\xi)
+ \pi_2^*d\te(\xi)dH_1(z)\eta
- dH_1(z)\xi \pi_2^*d\te(\eta) \\
& = (\om_W)_z(\xi,\eta),
\end{align*}
as desired.
\end{proof}

Now, we claim that the quotient projection $\pi_\mu: \JJ^{-1}(\mu) \to
\JJ^{-1}(\mu)/S^1$ is given by $\pi_1 \circ \psi^{-1}|_{\JJ^{-1}(\mu)}$, where
$\pi_1: P \to N$ is the projection onto the first factor. In effect, note that,
by the previous lemma, $\psi$ sends the orbits of $X_\JJ$ to the orbits of
$X_\Phi$. On the other hand, $\pi_1$ is the quotient projection at
$\Phi^{-1}(\mu) = \psi^{-1}(\JJ^{-1}(\mu))$ with respect to the trivial
circle bundle given by the orbits of $X_{\Phi}$.

To show that $\sigma_\mu$ is equal to $\tau$, note that
\begin{align*}
\pi_\mu^*\tau & = (\psi^{-1})*i_{\Phi^{-1}(\mu)}^*\pi_1^*\tau \\
& = (\psi^{-1})^*i_{\Phi^{-1}(\mu)}^*\om_W \\
& = i_{\psi(\Phi^{-1}(\mu))}^*\om_W \\
& = i_{\JJ^{-1}(\mu)}\om_W,
\end{align*}
where the third equality follows from the fact that $\psi^*\om_W=\om_W$.
\end{proof}

\subsection{Construction of a $\rho$-invariant pre-admissible Hamiltonian $\HH$ on $N \times T^*S^1$ from $H$. }
Consider the smooth map $\ppi: P \to N/\Z_n$ defined by
$$ \ppi|_{J^{-1}(\mu)} = \pi_\mu, $$
where $\pi_\mu$ is the quotient projection. Now, let $\bar H:= H\circ\tau^{-1}$ be the induced Hamiltonian on $N/\Z_n$. We will construct from $\bar H$ a $\rho$-invariant Hamiltonian defined on $P$. Fix a sufficiently small constant $\delta>0$ and define this new Hamiltonian by
$$ \HH(z) = (\bar H + \al(J(z))(m(\bar H)-\bar H))(\ppi(z)), $$
where $\al: \R \to \R$ is a $C^\infty$ function such that $0 \leq \al \leq 1$,
$\al(\mu) = 1\ \ \forall \mu \notin (\delta,1-\delta)$, $\al(\mu)= 0\ \ \forall
\mu \in (\frac{1}{2}-\delta,\frac{1}{2}+\delta)$ and $|\al^\pr(\mu)| \leq
2+2\delta$ for every $\mu \in [0,1]$. Thus, $\HH|_{J^{-1}(\mu)}$ is the lift of 
$\bar H_\mu:= \bar H + \al(\mu)(m(\bar H)-\bar H) = (1-\al(\mu))\bar H + \al(\mu)m(\bar H)$ by the quotient
projection $\pi_\mu$.

\begin{figure}[h]
\begin{center}
\includegraphics[width=2in]{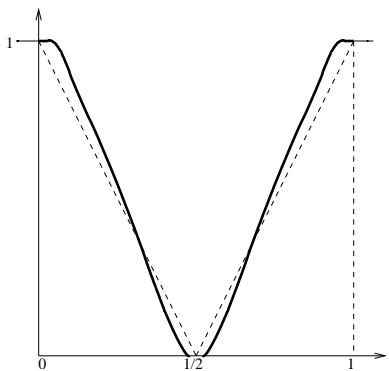}
\caption{\label{fctalpha} Graph of the function $\al$.}
\end{center}
\end{figure}

The first obvious property of $\HH$ is that $m(\HH) = m(\bar H) = m(H)$. In effect, $m(\HH)
= \sup_{\mu \in \R} m(\bar H+\al(\mu)(m(\bar H)-\bar H)) = m(\bar H)$, since $0 \leq \al \leq 1$
and $\al(\mu) = 1\ \ \forall \mu \notin (\delta,1-\delta)$.

It is easy to see that $\HH \in \H(\U,\om)$ where $\U$ is given by (see figure 3)
$$ \U = \bigcup_{0\leq\mu\leq 1} \pi_\mu^{-1}(\tau_n(U)).$$
In effect, it is clear that $0 \leq \HH \leq m(\HH)$. Moreover, $\HH|_{\V}
\equiv 0$ where $\V$ is the subset of $P$ given by
$\bigcup_{\frac{1}{2}-\delta\leq\mu\leq\frac{1}{2}+\delta} \pi_\mu^{-1}(\tau_n(V))$ and
$\tau_n(V) \subset N/\Z_n$ is the open set such that $\bar H|_{\tau_n(V)} \equiv 0$ because $\al(\mu) = 0\ \ \forall \mu \in
(\frac{1}{2}-\delta,\frac{1}{2}+\delta)$. Note that the interior of $\V$ is not
empty.

On the other hand, $\HH \equiv m(\HH)$ outside the compact set $\K =
\bigcup_{\delta\leq\mu\leq 1-\delta} \pi_\mu^{-1}(\tau_n(K)) \subset
\U\setminus\partial\U$, where $\tau_n(K) \subset \tau_n(U)\setminus\partial \tau_n(U)$ is the compact
set such that $\bar H|_{\tau_n(U)\setminus \tau_n(K)} \equiv m(\bar H)$.

Finally, note that if $H$ has no overtwisted critical points then $\HH$ also has no overtwisted critical points because $\bar H_\mu$ has no such critical points and the derivative of $X_\HH$ equals the identity in the direction tangent to the orbits of $\rho$.

\subsection{Existence of a non-trivial periodic orbit of $X_{\bar H}$ given by the reduction of a periodic orbit of $X_\HH$. }
Let us now prove that there exists a non-trivial periodic orbit $\bar\ga$ of $X_{\bar H}$ with period
less than $(2\pi/n)(\|H_1|_U\|+\sqrt{n})/m(H)$ and homotopy class $[\bar\ga] \in (\tau_n)_*G_\vr$. We will need the following theorem whose proof relies on results of Hofer-Viterbo \cite{HV} and McDuff-Slimowitz \cite{MSl} generalized for geometrically bounded symplectic manifolds by G. Lu \cite{Lu1,Lu2}.

\begin{theorem}
\label{HV}
Let $(M,\om)$ be a geometrically bounded symplectic manifold and $P:= M \times T^*S^1$
endowed with the symplectic form $\om_P:= \om \oplus \om^{S^1}_0$, where
$\om^{S^1}_0$ is the canonical symplectic form on $T^*S^1$. Let $U \subset P$
be an open subset with compact closure and $\tau_2: P \to T^*S^1$ the
projection onto the second factor. Then given a Hamiltonian $H \in
\H(U,\om_P)$ without overtwisted critical points, there exists a nonconstant periodic orbit $\gamma$ of $X_H$ with
period
$$ T < T_{max}(U,H) := \frac{\bigg|\int_{\tau_2(U)}\om^{S^1}_0\bigg|}{m(H)}. $$
Moreover,
the homotopy class $[\ga]$ of $\ga$ belongs to the subgroup $\pi_1(S^1)
\subset \pi_1(P)$.
\end{theorem}

\begin{proof}
Since $\ov U$ is compact, there exists a positive constant $a > 0$ such that
$\tau_2(U) \subset S^1 \times [-a/2,a/2] \subset T^*S^1$.  Let $\om_0$ be the
canonical symplectic form on $\R^2$ and consider the symplectomorphism $\phi:
S^1 \times [-a,a] \to A := \{(x,y) \in \R^2; a \leq x^2+y^2 \leq 5a \}$ given
by
$$ \phi(\theta,r) = (\sqrt{3a+2r}\sin\theta,\sqrt{3a+2r}\cos\theta). $$
Note that $Ker(Id,\phi)_*= \pi_1(S^1)$, where $(Id,\phi)_*: \pi_1(T^*M
\times S^1 \times [-a,a]) \simeq \pi_1(T^*P) \to \pi_1(T^*M \times \R^2) \simeq
\pi_1(T^*M)$ is the homomorphism induced on the fundamental group by the
transformation given by the identity and $\phi$ in the first and second factors
respectively.

Now, notice that it is sufficient to prove the theorem for an open subset $U^\prime$ such that $U \subset U^\prime$ and
$$ \bigg|\int_{\tau_2(U^\prime)}\om^{S^1}_0 - \int_{\tau_2(U)}\om^{S^1}_0\bigg| < \epsilon, $$
for $\epsilon>0$ arbitrarily small.

Thus, we can suppose, without loss of generality, that $\phi(\tau_2(\ov U)) \subset \R^2$ is a connected two-dimensional compact submanifold with boundary, such that there exists an open disk of radius $R$ with $L$ distinct points $y_j \in B^2(R)$ ($0 \leq L < \infty$) and an orientation preserving diffeomorphism
$$ \psi: \phi(\tau_2(U)) \to B^2(R)\setminus\{y_1,...,y_L\} $$
such that
$$ \bigg|\int_{\phi(\tau_2(U))}\om_0\bigg| = 
\bigg|\int_{B^2(R)\setminus\{y_1,...,y_L\}}\om_0\bigg|
= \pi R^2.$$

From a theorem of Dacorogna and Moser (see Lemma 2.2 in \cite{Si}) this $\psi$
can be required to be symplectic. Thus, consider the Hamiltonian $\bar H: M
\times B^2(R) \to \R$ given by
\begin{equation*}
\bar H(x,y)= 
\begin{cases}
m(H) \text{ if } y=y_j \text{ for some } j = 1,...,L \\
H(x,\phi^{-1}\psi^{-1}y) \text{ otherwise}
\end{cases}
\end{equation*}
which is obviously $C^\infty$, since $H|_{U\setminus K}=m(H)$, where $K \subset
U$ is a compact subset such that $K \subset U\setminus\partial U$.

By Proposition 1.6 of \cite{MSl} (extended for geometrically bounded symplectic manifolds in \cite{Lu1}), $X_{\bar H}$ has a contractible periodic orbit $\bar\ga$ with period
$$ T < \frac{\pi R^2}{m(\bar H)} =
\frac{\big|\int_{\phi(\tau_2(U))}\om_0\big|}{m(\bar H)} =
\frac{\big|\int_{\tau_2(U)}\om^{S^1}_0\big|}{m(H)}. $$

Finally, note that, by the remark above, the periodic orbit $\ga$ of $X_H$ given
by $\ga = (Id,\psi\circ\phi)^{-1}\bar\ga$ has homotopy class $[\ga]$ contained
in $\pi_1(S^1) \subset \pi_1(P)$, since $\bar\ga$ is contractible.
\end{proof}

\begin{remark}
If $m(M,\om) \geq \big|\int_{\tau_2(U)}\om^{S^1}_0\big|$, the hypothesis that $H$ has no overtwisted critical points is not necessary, see \cite{HV,Lu1,Lu2}.
\end{remark}

Now, notice that 
$$ \bigg|\int_{\tau_2(\U)}\om_0^{S^1}\bigg| = 2\pi\bigg(\frac{\|H_1|_U\|}{n} + \frac{1}{\sqrt{n}}\bigg) = (2\pi/n)(\|H_1|_U\| + \sqrt{n}). $$
In effect, the image of $\U \cup J^{-1}(0)$ under the projection $N \times S^1 \times \R \to N \times \R$ is the graph of $-(1/n)H_1$ restricted to $U$
and the projection of $\U \cup J^{-1}(1)$ is the graph of $1/\sqrt{n}-(1/n)H_1|_U$:

\begin{figure}[h]
\label{UU}
\begin{center}
\psfrag{mu}{$\mu$}
\psfrag{UU}{$\U$}
\psfrag{grH1U}{$-\frac{1}{n}H_1|_U$}
\psfrag{Phi1}{$\Phi^{-1}(1/\sqrt{n})$}
\psfrag{Phi0}{$\Phi^{-1}(0)$}
\psfrag{J1}{$J^{-1}(1)$}
\psfrag{J0}{$J^{-1}(0)$}
\includegraphics[width=2in]{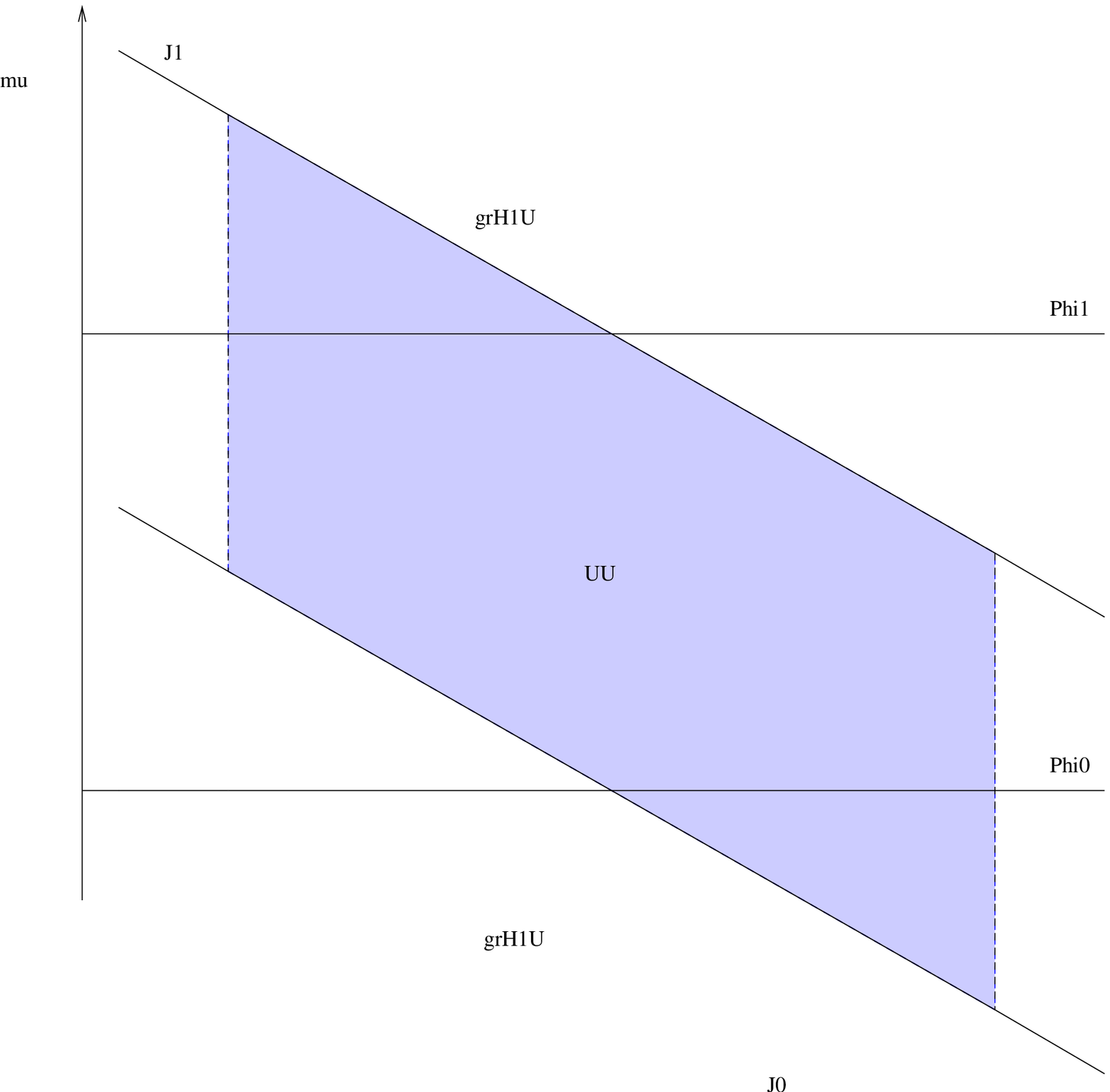}
\caption{\label{subsetU} The subset $\U$.}
\end{center}
\end{figure}

\vskip .3cm
\noindent {\bf Reduced dynamics of $X_\HH$:}
\vskip .2cm

Note that $\bar H_\mu = (1-\al(\mu))\bar H + \al(\mu)m(H)$ has the Hamiltonian vector
field with respect to $\om_n$ given by $X_{\bar H_\mu} = (1-\al(\mu))X_{\bar H}$. Thus, the
reduced dynamics of $X_\HH$ at $J^{-1}(\mu)$ is a reparametrization of the
dynamics of $X_{\bar H}$. Since $0 \leq \al(\mu) \leq 1$, the nonconstant periodic
orbits of $X_{\bar H_\mu}$ have period greater than or equal to those of $X_{\bar H}$.

Consequently, by Section \ref{step1}, it is sufficient to show
that the periodic orbit of $\HH$ given by Theorem \ref{HV} is not
tangent to an orbit of $\vr$ (such that it is projected
onto a singularity of $X_{\bar H}$).

\vskip .3cm
\noindent {\bf Non-triviality of the projected orbit:}
\vskip .2cm

Suppose that there exists a periodic orbit $\gga$ of $X_{\HH}$ tangent to an orbit of $\rho$. We will show that the period of $\gga$ is strictly greater than $T_{\max}(\U,\HH)$ (see Theorem \ref{HV}).

Firstly, note that
$$ T_{\max}(\U,\HH) = \frac{\big|\int_{\tau_2(\U)}\om_0^{S^1}\big|}{m(\HH)} =
\frac{(2\pi/n)(\|H_1|_U\|+\sqrt{n})}{m(H)}. $$

Now, let $W$ be an almost complex structure on $P$ compatible with $\om_P$ and
$\lg\cdot,\cdot \rg = \om_P(W\cdot,\cdot)$ be the corresponding Riemannian
metric. We have that along $\ga$,
$$ X_\HH = \bigg\lg X_\HH,\frac{X_J}{\|X_J\|^2}\bigg\rg X_J. $$
Consequently, since the period of the orbits of $\rho$ is equal to $\sqrt{n}$, the period of $\ga$ is given by
$$ T = \frac{\sqrt{n}}{\big|\big\lg X_\HH,\frac{X_J}{\|X_J\|^2}\big\rg\big|}. $$
But, by the definition of $n$,
$$ T = \frac{\sqrt{n}}{\big|\big\lg X_\HH,\frac{X_J}{\|X_J\|^2}\big\rg\big|} >
\frac{(2\pi/n)(2+2\delta)(\|H_1|_U\|+\sqrt{n})}{\big|\big\lg X_\HH,\frac{X_J}{\|X_J\|^2}\big\rg\big|}, $$
for $\delta>0$ sufficiently small.
Thus, it is sufficient to prove that
$$ \bigg|\bigg\lg X_\HH,\frac{X_J}{\|X_J\|^2}\bigg\rg\bigg| \leq
(2+2\delta)m(H), $$
such that
$$ T > \frac{(2\pi/n)(\|H_1|_U\| + \sqrt{n})}{m(H)}. $$
In effect, note that
\begin{align*}
\bigg|\bigg\lg X_\HH,\frac{X_J}{\|X_J\|^2} \bigg\rg\bigg|
& = \bigg|\bigg\lg W\nabla\HH,\frac{W\nabla J}{\|W\nabla J\|^2} \bigg\rg\bigg| \\
& = \bigg|\bigg\lg \nabla\HH,\frac{\nabla J}{\|\nabla J\|^2} \bigg\rg\bigg| \\
& = \bigg|d\HH\bigg(\frac{\nabla J}{\|\nabla J\|^2}\bigg)\bigg|,
\end{align*}
where in the second equation we used the fact that $W$ defines an
isometry. On the other hand, we have that for every $\xi \in T_zP$,
\begin{align*}
d\HH(z)\xi & =
d\bar H(\ppi(z))d\ppi(z)\xi
+ m(H)\al^\pr(J(z))dJ(z)\xi
 - \bar H(\ppi(z))\al^\pr(J(z))dJ(z)\xi\\
& \ \ - \al(J(z))d\bar H(\ppi(z))d\ppi(z)\xi\\
& = (m(H) - \bar H(\ppi(z)))\al^\pr(J(z))dJ(z)\xi
+ (1-\al(J(z)))d\bar H(\ppi(z))d\ppi(z)\xi.
\end{align*}
But note that, since $\gga$ is tangent to a fiber, it is projected onto a
singularity of $X_{\bar H}$, that is, $d\bar H(\ppi(z)) = 0$ for every $z \in \gga$.
Hence,
\begin{align*}
\bigg|\bigg\lg X_\HH,\frac{X_J}{\|X_J\|^2} \bigg\rg\bigg|
& = \bigg|d\HH\bigg(\frac{\nabla J}{\|\nabla J\|^2}\bigg)\bigg| \\
& = (m(H) - \bar H(\ppi(z)))|\al^\pr(J(z))| \\
& \leq (2+2\delta)m(H),
\end{align*}
as desired.

Consequently, by Section \ref{step1}, we conclude by symplectic reduction that $\bar H$ has a non-trivial periodic orbit $\bar\ga$ of period less than or equal to
$$ \frac{(2\pi/n)(\|H_1|_U\|+\sqrt{n})}{m(H)}. $$
But, since $\tau: N \to N/\Z_n$ is a symplectic covering with $n$ sheets, we have that the lift $\ga$ of $\bar\ga$ is a periodic orbit of $X_H$ with period less than or equal to
$$ \frac{2\pi(\|H_1|_U\|+\sqrt{n})}{m(H)}. $$

Finally, note that the periodic orbit $\gga$ of $X_\HH$ given by
Theorem \ref{HV} satisfies $[\gga] \in \pi_1(S^1)$. Hence,
the corresponding reduced periodic orbit $\bar\ga$ satisfies $[\bar\ga] \in (\tau_n)_*G_\vr$ since
$(\pi_\mu)_*\pi_1(S^1) = (\tau_n)_*G_\vr$. Thus, since the periodic orbit $\ga$ of $X_H$ is given by the lift of $\bar\ga$ by $\tau_n$, we have that $[\ga] \in G_\vr$.

\section{Proof of Theorem \ref{thm2}}
\label{proofthm2}

It is sufficient to prove Theorem \ref{thm2} for $n=2$, since proceeding inductively we can apply the theorem below for $N/\Z_n$ and conclude the result for any $n = 2^m$.

\begin{theorem}
Let $M$ be a symplectic manifold endowed with a $\Z_2$ free symplectic action, that is, a symplectomorphism $\psi: M \to M$ such that
$\psi^2=Id$ and $\psi(x)\neq x$ for every $x \in M$.  Then, given an admissible Hamiltonian $H$ on $M$ (resp. $H \in \Ht_a(M,\om)$), there exists an admissible Hamiltonian $H^\prime$ (resp. $H^\prime \in \Ht_a(M,\om)$) such that
\begin{itemize}
\item $H^\prime$ is $\Z_2$-invariant;
\item $\text{supp }X_{H^\prime} = \text{supp }X_{H+H\circ\psi}$;
\item $\|H^\prime\| \geq (1/2)\|H\|$.
\end{itemize}
Moreover, if $\psi$ is isotopic to the identity and $H$ is $G$-admissible (resp. $H \in \Ht_a^G(M,\om)$) then $H^\prime$ is also $G$-admissible (resp. $H^\prime \in \Ht_a^G(M,\om)$), for any subgroup $G \subset \pi_1(M)$.
\end{theorem}

\begin{proof}
Define $\HH = (1/2)(H + H\circ\psi)$. The first step in the proof is the following key proposition:

\begin{proposition}
Given a periodic orbit $\gamma: [0,T] \to M$ of $X_\HH$ then $X_H$ and $X_\Hp$
are colinear along $\gamma$.
\end{proposition}

\begin{proof}
Initially, note that $X_H$ and $X_\Hp$ commute, since
\begin{align*}
\{H,\Hp\} & = \om(X_H,X_\Hp) \\
& = \om(X_H,\psi_*X_H) \\
& = \om(\psi_*X_H,\psi_*\psi_*X_H) \\
& = \om(\psi_*X_H,X_H) = 0,
\end{align*}
where in the second and third equalities we used the fact that $\psi$ is
symplectic and the last equality follows by the antisymmetry of $\om$.

\begin{lemma}
\label{perorb}
If $X_H$ and $X_\Hp$ are not colinear along $\gamma$, then $\ga$ is not
isolated. In fact, there exists a 1-parameter deformation of $\ga$ by periodic
orbits of $X_\HH$ with the same period and energy, that is, there exists a family of
geometrically distinct periodic orbits $\gamma_s: [0,T] \to M$ for $s \in
(-\ep,\ep)$ for some $\ep > 0$ such that $\ga_0=\ga$ and $\HH(\ga_s(t)) = \HH(\ga_0(t))$ for every $s \in (-\ep,\ep)$.
\end{lemma}

\begin{proof}
Suppose that there exists $t \in [0,T]$ such that $X_H(\gamma(t))$ and
$X_\Hp(\gamma(t))$ are linearly independent. Then, $X_H(\gamma(t))$ and
$X_\Hp(\gamma(t))$ are not colinear to $X_\HH(\gamma(t))$. Since
$[X_H,X_\HH]=0$ we have that $\phi_s \circ \gamma: [0,T] \to M$ is a periodic
orbit for every $s \in \R$, where $\phi_s$ is the flow of $X_H$. These periodic
orbits are geometrically distinct for every $s$ sufficiently small because
$X_H(\gamma(t))$ and $X_\HH(\gamma(t))$ are linearly independent.
\end{proof}

Consider now the following linear application:
$$ \Phi: C^\infty(M) \to C^\infty(M/\Z_2) $$
given by $\Phi(H)(x) = (1/2)(H+\Hp)(\pi^{-1}(x))$, where $\pi: M \to M/\Z_2$ is
the quotient projection. Note that we can identify $C^\infty(M/\Z_2)$ with the
subspace $S \subset C^\infty(M)$  of $\psi$-invariant smooth functions and,
with respect to this identification, $\Phi$ is just the projection.

Now, suppose, by contradiction, that $X_H$ and $X_\Hp$ are not colinear along
$\ga$. By the previous lemma, we have that $\pi\circ\ga$ is not an isolated
periodic orbit of $X_{\Phi(H)}$. In fact, there exists a 1-parameter
deformation of $\pi\circ\ga$ given by periodic orbits of the same period.

\begin{lemma}[Lemma 19 of \cite{Rob}]
Given $\ep>0$, there exists $\FF \in C^\infty(M/\Z_2)$ such that $\|\FF -
\Phi(H)\|_{C^\infty} < \ep$ and  $\pi\circ\ga$ is an isolated periodic orbit of $X_\FF$ on the corresponding energy level.
\end{lemma}

Now, let $\FFF$ be the corresponding function in $S \subset C^\infty(M)$ given by the identification of $S$ with
$C^\infty(M/\Z_2)$. Write $\FFF = \HH + G$ such that $\|G\|_{C^\infty} <
\ep$ and define $F=H+G$. Note that obviously,
$$ \Phi(F) = \FF. $$

Since $\|F-H\|_{C^\infty}<\ep$ we have that $X_F$ and $X_{F\circ\psi}$ are linearly independent at $\ga(t)$. Consequently, by Lemma \ref{perorb}, there exists a deformation  of $\pi\circ\ga$ given by periodic orbits of $X_{\FF}$ with the same period and energy, a contradiction.
\end{proof}

\begin{proposition}
If $H$ has no overtwisted critical points then $\HH$ has no overtwisted critical points too.
\end{proposition}

\begin{proof}
Suppose, by contradiction, that $\HH$ has an overtwisted critical point $p$. By Corollary 3.5 of \cite{MSl}, there exists a sequence of Hamiltonians $\HH_n$ such that $\|\HH-\HH_n\|_{C^\infty} \overset{n\to\infty}\lra 0$ and $X_{\HH_n}$ has a non-trivial periodic orbit $\gga_n$ of period less than 1 converging to $p$. By the discussion in the proof of the previous proposition, there exists a sequence $H_n$ such that $\|H-H_n\|_{C^\infty} \overset{n\to\infty}\lra 0$ and either $X_{H_n}$ or $X_{H_n\circ\psi}$ has a non-trivial periodic orbit $\ga_n$ of period $<1$ given by a reparametrization of $\gga_n$. Thus, $p$ is also a critical point for both $H$ and $H\circ\psi$.

Now, note that, since $X_H$ and $X_\Hp$ commute, $DX_H(p)$ and $DX_\Hp(p)$ also commute. Hence, if $(1/2)(DX_H+DX_\Hp)(p)$ has an eigenvalue $\pm i\lambda$ for $\lambda>2\pi$, the same holds either for $DX_H(p)$ or $DX_\Hp(p)$, a contradiction.
\end{proof}

Now, let $\delta>0$ be a sufficiently small constant and $f: [\min \HH,\max \HH]
\to \R$ a smooth function such that $|f^\prime| \leq 1+2\delta$, $f(x)=0$ for every $x \in [\min \HH, \min \HH + \delta]$ and $f(x)=\|\HH\|$ for every $x \in [\max \HH-\delta, \max \HH]$:

\begin{figure}[h]
\label{functionf}
\begin{center}
\psfrag{medH}{$\|\HH\|$}
\psfrag{minH}{$\min \HH$}
\psfrag{maxH}{$\max \HH$}
\includegraphics[width=2in]{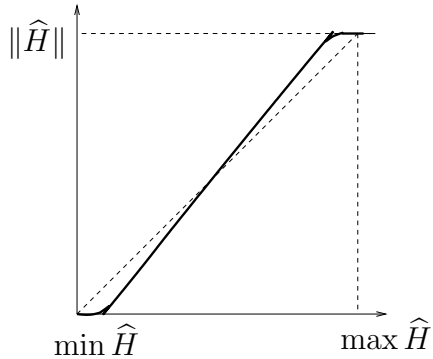}
\caption{\label{fctf} The function $f$.}
\end{center}
\end{figure}

Consider the Hamiltonian given by
$$ H^\prime = f\circ \HH. $$
Note that, since $\psi$ is proper and by the definition of $f$, $H^\prime$ is pre-admissible. Moreover, we
have that $\|H^\prime\| = \|\HH\| \geq (1/2)\|H\|$, because $\min \HH \leq
(1/2)\max H$ and $\max \HH = \max H$, because $\psi$ is proper.

Now, suppose that $X_{H^\prime}$ has a periodic orbit of period $T^\prime < 1$ and homotopy class in $G$.
Then, since $X_{H^\prime}(x) = f^\prime(\HH(x))X_\HH(x)$ and $|f^\prime| \leq
1+2\delta$ for $\delta$ arbitrarily small, we have that $X_\HH$ also has a periodic
orbit $\ga$ of period $T<1$ such that $[\ga] \in G$.

By the previous proposition, $X_H$ and $X_\Hp$ are colinear to $X_\HH$ along
$\ga$. Consequently, since $X_\HH = (1/2)(X_H + X_\Hp)$, we conclude that $\ga$
is also a periodic orbit of period less than 1 either for $X_H$ or $X_\Hp$.
But, how $H$ and $\Hp$ are both admissible (resp. $G$-admissible) because $\psi$ is symplectic (resp. because $\psi$ is isotopic to the identity), we have a contradiction.

\end{proof}

\end{document}